\documentclass[11pt]{amsart}

\usepackage{amsthm,amsmath,amssymb,graphics,graphicx,hyperref,epstopdf}

\title{Generalizing the Borel Property}
\author{Christopher A. Francisco}
\address{Department of Mathematics, Oklahoma State University,
401 Mathematical Sciences, Stillwater, OK 74078}
\email{chris@math.okstate.edu}
\urladdr{http://www.math.okstate.edu/$\sim$chris}

\author{Jeffrey Mermin}
\address{Department of Mathematics, Oklahoma State University,
401 Mathematical Sciences, Stillwater, OK 74078}
\email{mermin@math.okstate.edu}
\urladdr{http://www.math.okstate.edu/$\sim$mermin}

\author{Jay Schweig}
\address{Department of Mathematics, University of Kansas, 405 Snow Hall, Lawrence, KS 66045}
\email{jschweig@math.ku.edu}
\urladdr{http://www.math.ku.edu/$\sim$jschweig}

\newtheorem{theorem}{Theorem}[section]
\newtheorem{proposition}[theorem]{Proposition}
\newtheorem{corollary}[theorem]{Corollary}
\newtheorem{lemma}[theorem]{Lemma}

\theoremstyle{definition}
\newtheorem{definition}[theorem]{Definition}
\newtheorem{remark}[theorem]{Remark}
\newtheorem{example}[theorem]{Example}

\newtheorem{notation}[theorem]{Notation}
\newtheorem{question}[theorem]{Question}
\newtheorem{algorithm}[theorem]{Algorithm}

\DeclareMathOperator{\Q}{Q}

\DeclareMathOperator{\codim}{codim}
\DeclareMathOperator{\Ann}{Ann}
\DeclareMathOperator{\Ass}{Ass}

\DeclareMathOperator{\pd}{pd}
\DeclareMathOperator{\dist}{dist}
\DeclareMathOperator{\Beg}{beg}
\DeclareMathOperator{\END}{end}
\DeclareMathOperator{\supp}{supp}
\DeclareMathOperator{\Hilb}{Hilb}
\DeclareMathOperator{\Borel}{Borel}

\newcommand{\p}{\mathfrak{p}}
\newcommand{\q}{\mathfrak{q}}

\begin{document} 

\begin{abstract}
We introduce the notion of $Q$-Borel ideals:  ideals which are closed under the Borel moves arising from a poset $Q$.  We study decompositions and homological properties  of these ideals, and offer evidence that they interpolate between Borel ideals and arbitrary monomial ideals.
\end{abstract}

\maketitle

\section{Introduction} \label{s.introduction}

Borel-fixed ideals are a natural class of ideals to study in commutative algebra: Not only do they arise in a host of contexts, but they are typically easier to understand than arbitrary monomial ideals. When $k$ is a field of characteristic zero, an ideal $B \subseteq k[x_1, x_2, \ldots, x_n]$ is Borel-fixed if and only if it satisfies the following property: If $m$ is a monomial in $B$ and $x_j$ divides $m$, then $m \cdot \frac{x_i}{x_j} \in B$ for all $i < j$. That is, $B$ is closed under ``Borel moves.'' Over an arbitrary field, we call the ideals for which this combinatorial condition holds \emph{Borel ideals}, and our previous paper \cite{FMS} was devoted to studying Borel ideals via their \emph{Borel generators}.  

In this paper, we investigate how loosening the restrictive Borel condition affects the properties of a monomial ideal. That is, what happens if we require only that $m \cdot \frac{x_i}{x_j}$ remain in the ideal for some pairs $x_i,x_j$ rather than all pairs with $i<j$? We formalize this framework by introducing the notion of a monomial ideal being \emph{Borel with respect to $Q$}, where $Q$ is a poset on $\{x_1, \dots, x_n\}$. A monomial ideal $I$ is \emph{$Q$-Borel} if for any monomial $m \in I$ with $x_j$ dividing $m$ and any $x_i <_Q x_j$, we have $m \cdot \frac{x_i}{x_j} \in I$ (here $<_Q$ denotes the relation in the poset $Q$).

Two extremal cases are notable: When $Q$ is the chain $C:  x_1 <_Q x_2 <_Q \cdots <_Q  x_n$, the $C$-Borel ideals are precisely the usual Borel ideals. When $Q$ is the antichain, every monomial ideal is $Q$-Borel because there are no conditions to satisfy. This perspective is helpful because it allows us to use the Borel approach to study more general ideals satisfying a subset of the Borel conditions.

If a monomial ideal $I$ is $Q$-Borel, and $Q$ is similar to the chain $C$, then $I$ should behave similarly to a Borel ideal. 

In particular, if $Q$ is ``close'' to the chain and $I$ is $Q$-Borel, then $I$ should have minimal free resolution that resembles the Eliahou-Kervaire resolution, and the associated primes of $I$ should be similar to the associated primes of a Borel ideal. We are especially interested in understanding the \emph{principal} $Q$-Borel ideals, namely the ideals which have a single $Q$-Borel generator and whose ordinary monomial generating set arises from performing all possible sequences of $Q$-Borel moves on that generator.

Our work also has connections to matrix groups. Let $U_n(k)$ be the Borel subgroup of $GL_n(k)$, consisting of the set of invertible $n \times n$ upper-triangular matrices with entries in the field $k$.  That is, $U_n(k) = \{A \in GL_n(k): a_{i,j} = 0$ for $i > j\}$.  Let $D_n(k)$ denote the set of invertible diagonal $n \times n$ matrices with entries in $k$.  Note that $D_n(k) \subseteq U_n(k)$.  

If $A$ is any $n \times n$ matrix with entries in $k$, recall that $A$ acts on polynomials in $k[x_1, x_2, \ldots, x_n]$ by replacing each occurrence of $x_i$ with $\sum_{j = 1}^n a_{ij} x_j$.  The classical Borel ideals are exactly those that are fixed under the action of $U_n(k)$, whereas monomial ideals are exactly those fixed by $D_n(k)$ (see, for instance, \cite{MS}).  Now let $Q$ be a naturally labeled poset on $\{x_1, x_2, \ldots, x_n\}$, and define $M_Q(k)$ as follows:
\[
M_Q (k)= \{A \in GL_n(k): a_{i, j} = 0 \text{ whenever } x_i \nleq_Q x_j\}.
\]
Because $Q$ is naturally labeled, we have $D_n(k) \subseteq M_Q(k) \subseteq U_n(k)$.  Moreover, we have the following: 

\begin{proposition}
The $Q$-Borel ideals in $k[x_1, \ldots, x_n]$ are exactly those which are fixed under the action of $M_Q(k)$.  
\end{proposition}

A common theme throughout the paper will be that, when $Q$ is the chain or the antichain, the theory of $Q$-Borel ideals specializes to known results.  This heuristic applies here as well:  If $Q$ is the antichain, then $M_Q(k) = D_n(k)$, whereas $M_Q(k) = U_n(k)$ when $Q$ is the chain.  

The paper is organized as follows. We present some preliminary results about $Q$-Borel ideals in Section~\ref{s.preliminaries}. 
In Section \ref{s.primary}, we use M\"obius inversion and Hall's Marriage Theorem to study the relationship between products of monomial primes and monomial ideals which decompose as intersections of prime powers.  Our main results are Theorems \ref{primeproduct} and \ref{colontheorem}. Theorem \ref{primeproduct} explicitly decomposes a product of monomial primes into an intersection of prime powers. Theorem \ref{colontheorem} describes, in terms of colon ideals, all the ideals which can be decomposed as an intersection of powers of monomial primes. Section \ref{s.primaryprincipal} applies Theorem \ref{primeproduct} to produce irredundant primary decompositions of principal $Q$-Borel ideals. The results of this section were recently proved in a broader context by Herzog, Rauf, and Vladoiu \cite{HRV} in independent work. In Section~\ref{s.cm}, we compute the projective dimension and codimension of any principal $Q$-Borel ideal. As a corollary, we determine when principal $Q$-Borel ideals are Cohen-Macaulay, recovering part of a result of Herzog and Hibi \cite{herzoghibi}.  In Section \ref{s.interpolation}, we describe algorithms for computing free resolutions and irreducible decompositions of $Q$-Borel ideals.  These algorithms interpolate between familiar constructions for Borel ideals and for arbitrary monomial ideals.
In Section~\ref{s.y-borel}, we explicitly construct the minimal free resolutions of $Y$-Borel ideals, where $Y$ is a specific poset close to the chain.  This resolution is very similar to the Eliahou-Kervaire resolution, and provides evidence for the idea that, if $Q$ is a poset close to the chain, then $Q$-Borel ideals should behave similarly to Borel ideals.

\section{Preliminary results}\label{s.preliminaries}

Let $Q$ be a poset on $\{x_1, x_2, \ldots, x_n\}$ with partial order $<_Q$.  Recall that $Q$ is \emph{naturally labeled} if $x_i <_Q x_j$ implies $i < j$. All posets considered here will be naturally labeled.  Our use of poset terminology is standard; see \cite{ec1}.   \\

Throughout, let $S = k[x_1, x_2, \ldots, x_n]$ where $k$ is a field.  The central objects considered in this paper are $Q$-Borel ideals, defined as follows. 

\begin{definition}
Let $I \subseteq S$ be a monomial ideal, and let $Q$ be a poset on $\{x_1, x_2, \ldots, x_n\}$.  We say that $I$ is \emph{$Q$-Borel} if whenever $x_i <_Q x_j$ and $m\in I$ is a monomial divisible by $x_{j}$, then $m \cdot \frac{x_i}{x_j} \in I$.  In this case, we may also say that $Q$ \emph{stabilizes} $I$, or that $I$ is \emph{Borel with respect to $Q$}.    
\end{definition}

If $Q$ and $Q'$ are two posets, we say that $Q$ \emph{refines} $Q'$ if $x_i <_{Q'} x_j$ implies $x_i <_{Q} x_j$ for any $i$ and $j$.  

\begin{proposition}
Let $I \subseteq S$ be a monomial ideal.  Then there exists a poset $Q$ with the following properties.  
\begin{enumerate}
\item[(1):] $I$ is $Q$-Borel, and 
\item[(2):] If $Q'$ is a poset for which $I$ is $Q'$-Borel, then $Q$ refines $Q'$.  
\end{enumerate} 
We call the poset $Q$ the \emph{maximal stabilizing poset} of $I$.  
\end{proposition}

\begin{proof}
Suppose $Q_1$ and $Q_2$ both stabilize $I$.  Let $Q_3$ be the poset that is the transitive closure of the union of the relations in $Q_1$ and $Q_2$.  Then $Q_3$ stabilizes $I$ as well.  Continuing in this way, we eventually obtain a poset $Q$ with the desired properties.  
\end{proof}

Recall that an \emph{order ideal} of a poset $Q$ is a subset $A$ of its elements such that $y\in A$ and $x <_Q y$ implies $x\in A$. By abuse, if $A$ is an order ideal, we also let $A$ denote the prime ideal generated by the elements of $A$.

\begin{theorem}
For any poset $Q$, there exists a monomial ideal $I$ such that $Q$ is its maximal stabilizing poset.
\end{theorem}

\begin{proof}
Let $I=\prod_{A}A$ be the product of all the order ideals of $Q$. Then $Q$ stabilizes $I$. If $x_j \not \leq_Q x_i$, let $X$ be the set of order ideals containing $x_{j}$, $Y$ the set of order ideals containing $x_{i}$ but not $x_{j}$, and $Z$ the set of order ideals containing neither.  Put
$m = \prod_{A\in X}x_{j}\prod_{A\in Y}x_{i}\prod_{A\in Z}x_{A}$, where $x_{A}$ is any element of $A$.  Observe that $x_{i}$ divides $m$, since $Y$ contains the principal order ideal generated by $x_{i}$.  Furthermore, $m \cdot \frac{x_j}{x_i} \notin I$.
\end{proof}

\begin{definition}
Let $X \subseteq S$ be a set of monomials, and let $Q$ be a poset.  If $I$ is the smallest $Q$-Borel ideal containing $X$, we say that $I$ is \emph{generated as a $Q$-Borel ideal} by $X$ and write $I=\Q(X)$.  In this case, we call $X$ a $Q$-Borel generating set for $I$.  Every $Q$-Borel ideal $I$ has a unique minimal $Q$-Borel generating set; we call the monomials in this set the \emph{$\Q$-generators} of $I$.  
Of particular interest is the case when $X$ consists of a single monomial $m$.  We say in this case that $I$ is a \emph{principal $Q$-Borel ideal}, and write $I=\Q(m)$. 
\end{definition}

\begin{definition}
If $Q$ is a poset with $x_i \leq_Q x_j $ and $m$ is a monomial with $x_j|m$, we call the replacement of $m$ with $m \cdot \frac{x_i}{x_j}$ a \emph{Q-Borel move}. 
\end{definition}

Thus, we may alternately define a $Q$-Borel ideal as a monomial ideal closed under $Q$-Borel moves.

Note that if $I = \Q(m)$, then the (ordinary) monomial generators of $I$ are all of the same degree.  This is a special case of the following.

\begin{proposition}
Let $X$ be a set of monomials, all of the same degree.  Then the ordinary monomial generators of $I = \Q(X)$ are all of the same degree.  
\end{proposition}

\begin{proof}
The ordinary generators of $I$ are the monomials that can be obtained, via a sequence of $Q$-Borel moves, from monomials in $X$.  Because performing a $Q$-Borel move on a monomial cannot change its degree, the result follows. 
\end{proof}

\begin{proposition}\label{primefactor}
Let $I = Q(m)$ be a principal $Q$-Borel ideal.  Then $I$ can be factored as a product of monomial prime ideals.  Moreover, the primes in this factorization are all $Q$-Borel ideals, and can be determined from the exponent vector of $m$.
\end{proposition}
\begin{proof}
Write $m=\prod x_{i}^{e_{i}}$.  For a variable $x_{i}$, let $\frak{p}_{i}=(x_{j}:x_{j} \le_{Q}x_{i})$ be the order ideal of $Q$ generated by $x_{i}$.  We claim that 
\[
I=\prod_{x_i | m} \frak{p}_{i}^{e_{i}}.
\]
Indeed, suppose $\mu\in I$.  We may assume $\deg(\mu)=\deg(m)$.  Then $\mu$ can be obtained from $m$ by a sequence of $Q$-Borel moves, i.e., $\mu=m\cdot\prod\frac{w_{i}}{x_{i}}=\prod w_{i}$ for some collection $\{w_{i}\}$ of variables with $w_{i}\leq_{Q} x_{i}$ for all $i$, so that $w_{i}\in \frak{p}_{i}$ for all $i$, i.e., $\mu\in\prod \frak{p}_{i}$.

Conversely, suppose that $\mu\in \prod\frak{p}_{i}$.  Again we may assume $\deg(\mu)=\deg(m)$.  Write $\mu=\prod w_{i}$ for some collection $\{w_i\}$ of variables with $w_{i}\in \frak{p}_{i}$ for all $i$.  Then $\mu=\prod x_{i}\frac{w_{i}}{x_{i}} = m\prod\frac{w_{i}}{x_{i}}$ is obtained from $m$ by $Q$-Borel moves, meaning $\mu\in I$.
\end{proof}

\begin{definition}
Let $I$ be a monomial ideal with minimal generating set $G(I)$, where the monomials in $G(I)$ are all of the same degree.  Then $I$ is called \emph{polymatroidal} if the following exchange condition is satisfied:  For two monomials $m, m' \in G(I)$ and any variable $x_i$ appearing to a greater power in $m$ than in $m'$, there exists a variable $x_j$, appearing to a greater power in $m'$ than in $m$, such that $m \cdot \frac{x_j}{x_i} \in G(I)$.  
\end{definition}

\begin{proposition}\label{polym}
If $I$ is a principal $Q$-Borel ideal, then $I$ is polymatroidal.  
\end{proposition}

\begin{proof}
Each of the ideals $\frak{p}_i$ in Proposition \ref{primefactor} is prime, hence polymatroidal.  Theorem 5.3 of \cite{concaherzog}, which states that the product of polymatroidal ideals is polymatroidal, completes the proof. 
\end{proof}

\begin{definition}\label{linearquotients}
Suppose $I=(m_1,\dots,m_r) \subset S$, and the $m_i$ are monomials such that $\deg m_1 \le \cdots \le \deg m_r$. We say that $I$ has \emph{linear quotients} if, for each $1 < i \le r$, the ideal quotient $(m_1, \dots, m_{i-1}):m_i$ is generated by a subset of $x_1, \dots, x_n$.  If an ideal generated in a single degree has linear quotients, it has a linear resolution (see, for instance, \cite[Lemma 4.1]{concaherzog}).
\end{definition}

\begin{corollary}\label{linearresolution}
Every principal $Q$-Borel ideal has linear quotients and hence a linear resolution. 
\end{corollary}

\begin{proof}
In \cite[Lemma 1.3]{herzogtakayama}, Herzog and Takayama prove that polymatroidal ideals have linear quotients. Thus polymatroidal ideals generated in a single degree have a linear resolution.
\end{proof}

\section{Prime power decompositions} \label{s.primary}

We begin by proving a useful theorem for ideals that are products of monomial primes.  Note that Theorem~\ref{primeproduct} immediately produces a primary decomposition of such an ideal.  

\begin{theorem}\label{primeproduct}
Let $I \subseteq S$ be a product of monomial primes.  Write
\[
I = \prod_{\p} \p^{e_{\p}}, 
\]
where $\p$ ranges over the set of all monomial primes, and we allow any of the above exponents to be zero.  Then 
\[
I=\bigcap_{\p} \p^{a_{\p}},
\]
where each $\displaystyle{a_{\p} = \sum_{\q \subseteq \p} e_{\q}}$.  Furthermore, we have 
\[
e_{\p} = \sum_{\q \subseteq \p} (-1)^{|\p| - |\q|} a_{\q}.
\]
\end{theorem}

Before proving Theorem \ref{primeproduct}, we need a lemma which will use Hall's Marriage Theorem, a standard result in enumerative combinatorics (see, for instance, \cite{ec1}).  

Let $G$ be a bipartite graph on vertex set $X \sqcup Y$, so that each edge of $G$ contains one vertex in $X$ and one in $Y$.  For $A \subseteq X$, let $N(A)$ denote the set of neighbors of vertices in $A$.  That is, 
\[
N(A) = \{y \in Y : (x, y) \text{ is an edge of } G \text{ for some } x\in A\}.
\]
Recall that a \emph{perfect matching} of $G$ is a set of $|X|$ vertex-disjoint edges of $G$.

\begin{theorem}[Hall's Marriage Theorem]\label{marriage}
Let $G$ be as above.  Then $G$ has a perfect matching if and only if $|A| \leq |N(A)|$ for all $A \subseteq X$.
\end{theorem}

\begin{lemma}\label{monomialmarriage}
Let $I = \prod_{\p} \p^{e_\p}$, where each $e_{\p} \geq 0$.  Let $E = \sum e_{\p}$, and let $m$ be a monomial of degree $E$.  Then $m \in I$ if and only if, for each $\p$, the number of variables of $m$ in $\frak{p}$ is at least $\sum_{\q \subseteq \p} e_{\q}$. 
\end{lemma}

\begin{proof}  The ``only if'' direction is clear.  For the ``if'' direction let $m = \prod x_i^{f_i}$, and define a bipartite graph $G$ as follows.  Let $X$ be a set of $E$ vertices, $f_1$ of which are labeled $x_1$, $f_2$ of which are labeled $x_2$, and so on.  Similarly, let $Y$ be a set of $E$ vertices, $e_{\p}$ of which are labeled $\p$ for each $\p$. Now connect each vertex labeled $x_i$ to each vertex labeled $\p$ whenever $x_i \in \frak{p}$.  

We wish to apply Theorem \ref{marriage} to the graph $G$.  Indeed, if we can show that $|A| \leq |N(A)|$ for all $A \subseteq X$, the perfect matching guaranteed by the Marriage Theorem will give $m \in I$. 

Let $A \subseteq X$.  Then we can expand $A$ to include all vertices in $X$ that share labels with vertices in $A$ (this does not change the size of $N(A)$).  Thus $A =  \{x_i : x_i \notin \frak{p}\}$ for some $\p$, and 
\[
|A| = E - \sum_{x_i \in \frak{p}} f_i \leq E - \sum_{\frak{q} \subseteq \frak{p}} e_{\q}.
\]  
Similarly, we have $N(A) = \{\frak{q} : x_k \in \frak{q}$ for some $x_k \notin \frak{p}\} = \{\frak{q} : \frak{q} \nsubseteq \frak{p}\}$, and so
\[
|N(A)| = \sum_{\frak{q} \nsubseteq \frak{p}} e_{\q} = E - \sum_{\frak{q} \subseteq \frak{p}} e_{\q},
\]
so $|A| \leq |N(A)|$, completing the proof.  
\end{proof}

\begin{proof}[Proof of Theorem \ref{primeproduct}]
Let $m \in \prod \p^{e_{\p}}$, and assume $m$ is of degree $E$.  Then we can write $m = \prod m_{\p}$, where $\deg(m_{\p}) = e_{\p}$ for all $\p$ and each $m_{\p}$ has support in $\frak{p}$.  Now fix $\p$, and let $M=\prod_{\q \subseteq \p} m_{\q}$.  We have $m = M\prod_{\q \nsubseteq \p} m_{\q}$ and $\deg(M) = \sum_{\q \subseteq \p} e_{\q} = a_{\p}$, so $M \in \frak{p}^{a_{\p}}$, so $m \in \frak{p}^{a_{\p}}$.  Since the choice of $\p$ was arbitrary, we have $I \subseteq \bigcap \p^{a_{\p}}$. 

Now assume that $m \notin I$, and note we can assume that $\deg(m) = E$.  Then by Lemma \ref{monomialmarriage} there must be some $\frak{p}$ such that $m$ contains strictly fewer than $\sum_{\q \subseteq \p}e_{\q}$ variables from $\frak{p}$.  But then $m \notin \frak{p}^{a_{\p}}$, and so $m \notin \bigcap \p^{a_{\p}}$.  

For the final claim, apply the Principle of Inclusion-Exclusion.\end{proof}


\begin{corollary}
Let $I=\prod \frak{p}^{e_{\p}}$ be a product of monomial primes.
Then $\frak{q}\in \Ass(S/I)$ whenever $e_{\q}\geq 1$.
\end{corollary}
\begin{proof}
It suffices to show that the primary ideal $\q^{a_{\q}}$ is not
redundant in the decomposition $I=\bigcap \p^{a_{\p}}$.  Localizing, we
may assume that $\p$ is the homogeneous maximal ideal $(x_{1},\dots,
x_{n})$.  Now let 
$I'=\q^{e_{\q}-1}\displaystyle\prod_{\p\neq \q}\p^{e_{\p}}$
and apply Theorem \ref{primeproduct} to $I'$.  This yields
$I'=\q^{a_{q}-1}\cap \displaystyle\bigcap_{\p\neq \q}\p^{e_{\p}}$.
Since $I=\q^{a_{q}}\cap \displaystyle\bigcap_{\p\neq \q}\p^{e_{\p}}$
and $I\neq I'$, it follows that $\q^{a_{\q}}$ is not
redundant.
\end{proof}

For the following two corollaries we use M\"obius inversion (see \cite{ec1}). 

\begin{corollary}\label{anyposetprimarydecomp}
Let $I=\prod\frak{p}^{e_{\p}}$ be a product of monomial primes, and
let $\Lambda$ be any collection of monomial primes containing
$\Ass(S/I)$, ordered by inclusion.  Let $\mu$ be the M\"obius function
on $\Lambda$.  Then $I=\displaystyle\bigcap_{\frak{p}\in \Lambda}
\frak{p}^{a_{\p}}$, where
$a_{\p}=\displaystyle\sum_{\q\subseteq\frak{p}}e_{\q}$.  Furthermore, we can
recover the factorization of $I$ from this decomposition, via the formula
$e_{\p}=\displaystyle\sum_{\q\leq_{_{\Lambda}}\frak{p}}
\mu(\q, \frak{p})a_{\q}$.
\end{corollary}
\begin{proof}
We have $I=\displaystyle\bigcap_{\frak{p}\in
  \Lambda}\frak{p}^{a_{\p}}\cap \bigcap_{\frak{p}\not\in
  \Lambda}\frak{p}^{a_{\p}}$.  Since $\Lambda\supseteq \Ass(S/I)$, the
second factor is redundant and may be omitted, leaving the first
formula.  M\"obius inversion yields the second formula.

\end{proof}

\begin{proposition}\label{sumass}
Let $I=\prod \p^{e_{\p}}$ and $\q \in \Ass(S/I)$. Then $\q = \sum_{\p \in T} \p$, where $T$ is some set satisfying $\p \in T$ implies $e_{\p} \gneqq 0$.
\end{proposition}

\begin{proof}
Let $\displaystyle{\q' = \sum_{\stackrel{\p \subseteq \q}{e_{\p} \gneqq 0}} \p}$. Observe $\q' \subseteq \q$ and $a_{\q'} = a_{\q}$. Since $\q \in \Ass(S/I)$, we must have $\q = \q'$.
\end{proof}

\begin{proposition}
Suppose $\Lambda$ is closed under the taking of sums and that $I = \bigcap_{\p \in \Lambda} \p^{a_{\p}}$. Set $e_{\p}=\displaystyle\sum_{\frak{q} <_{\Lambda} \frak{p}} \mu(\frak{q},\frak{p})a_{\frak{q}}$. Suppose all $e_{\p} \ge 0$. Then $I = \prod_{\p \in \Lambda} \p^{e_{\p}}$.
\end{proposition}

\begin{proof}
By the previous proposition, $\Lambda$ contains $\Ass(S/I)$.  Apply Corollary \ref{anyposetprimarydecomp} to the ideal $\prod_{\p \in \Lambda} \p^{e_{\p}}$.
\end{proof}

\begin{theorem}\label{colontheorem}
Let $\Lambda$ be a set of monomial primes which is closed under the
taking of sums, and let $I$ have primary
decomposition 
\[
I=\bigcap_{\frak{p}\in \Lambda} \frak{p}^{a_{\frak{p}}}.
\]
Let $\mu$ be the M\"obius function on $\Lambda$, and set
$e_{\p}=\displaystyle\sum_{\frak{q}\subseteq\frak{p}}
\mu(\frak{q},\frak{p})a_{\frak{q}}$ for all $\p\in\Lambda$.
Then 
\[
I=\prod \p^{e_{\p}} = \left(\prod_{e_{\p}>0}\p^{e_{\p}} :
\prod_{e_{\p}<0} \p^{-e_{\p}}\right).
\]
\end{theorem}

\begin{proof}
Let $J=\prod_{e_{\p}>0}\p^{e_{\p}}$ and $K=\prod_{e_{\p}<0}
\p^{-e_{\p}}$. By Proposition \ref{sumass} and Corollary \ref{anyposetprimarydecomp}, we may write $J=\bigcap_{\p \in \Lambda} \p^{b_{\p}}$ and $K=\bigcap_{\p \in \Lambda} \p^{c_{\p}}$. Observe that $b_{\p} = a_{\p} + c_{\p}$ for all $\p$. Suppose $m \in (J:K)$. Fix $\p$. Define $\mu = \prod_{\q: e_{\q} \lneqq 0} x_{\q}^{-e_{\q}}$, where $x_{\q} \in \q$ is chosen so that $x_{\q} \notin \p$ whenever  $\q \not \subseteq \p$. Observe that $\mu \in K$ and $\mu \in \p^{c_{\p}} \smallsetminus \p^{c_{\p}+1}$. Because $m \mu \in J$, we have $m \mu \in \p^{b_{\p}}$, so $m \in \p^{b_{\p}-c_{\p}} = \p^{a_{\p}}$. Therefore $(J:K) \subseteq \p^{a_{\p}}$ for each $\p$, so $(J:K) \subseteq I$.

Now suppose $m \in I$. For all $\p$, we have $m \in \p^{a_{\p}} = \p^{b_{\p}-c_{\p}}$. Let $\mu \in K$. Then $m\mu \in \p^{b_{p}-c_{\p}+c_{\p}}=\p^{b_{\p}}$ for all $\p$. Consequently, $m \mu \in J$.
\end{proof}

\begin{corollary}
Let $\Lambda$ be a set of monomial primes ordered by inclusion, and let $I$ have primary
decomposition 
\[
I=\bigcap_{\frak{p}\in \Lambda} \frak{p}^{a_{\frak{p}}}.
\]
Let $\mu$ be the M\"obius function on $\Lambda$, and set
$e_{\p}=\displaystyle\sum_{\frak{q} <_{\Lambda} \frak{p}}
\mu(\frak{q},\frak{p})a_{\frak{q}}$ for all $\p\in\Lambda$.
Suppose $\Lambda'\supset \Lambda$ is another set of monomial primes (again ordered by inclusion),
and let $e_{\p}=0$ for all $\p\in\Lambda' \smallsetminus \Lambda$.  Then, if
$b_{\p}=\displaystyle\sum_{\q\leq_{_{\Lambda'}}\p}e_{\p}$, we have 
\[
I=\bigcap_{\frak{p}\in \Lambda'} \frak{p}^{b_{\frak{p}}}.
\]
\end{corollary}
\begin{proof}
Apply M\"obius inversion twice.
\end{proof}

\begin{question}
Let $\Lambda$ be a set of monomial primes which is closed under the
taking of sums, and let $I$ have primary
decomposition 
\[
I=\bigcap_{\frak{p}\in \Lambda} \frak{p}^{a_{\frak{p}}}.
\]
Let $\mu$ be the M\"obius function on $\Lambda$, and set
$e_{\p}=\displaystyle\sum_{\frak{q}\subseteq\frak{p}}
\mu(\frak{q},\frak{p})a_{\frak{q}}$ for all $\p\in\Lambda$.

Is there is an easily-checked condition on $\Lambda$ and the exponents $e_{\p}$ and $a_{\p}$ that determines which primes are associated?
\end{question}

\section{Primary decompositions of principal $Q$-Borel ideals} \label{s.primaryprincipal}

In this section, we apply Theorem \ref{primeproduct} and Proposition \ref{primefactor} to describe an irredundant primary decomposition of a principal $Q$-Borel ideal. Let $Q$ be a naturally labeled poset.  We say an order ideal of $Q$ is \emph{connected} if it cannot be written as the disjoint union of two order ideals.  

\begin{proposition}\label{qborelprimes}
Let $I$ be a $Q$-Borel ideal.  If $\frak{p}$ is an associated prime of $S/I$, then
\[
\frak{p} = (x_i : x_i \in A)
\]
for some order ideal $A$ of $Q$.  
\end{proposition}

\begin{proof}
Let $x_j \in \frak{p}$. Then there exists a monomial $m$ such that $m \notin I$, but $x_jm \in I$. If $x_i \le_Q x_j$, then by a $Q$-Borel move, $x_j m \cdot \frac{x_i}{x_j} = x_i m \in I$, so $x_i \in \frak{p}$ as well.
\end{proof}

\begin{corollary}\emph{\cite[Corollary 2]{bayerstillman}}\label{bayerstillman}
If $\frak{p}$ is an associated prime of a Borel ideal $S/I$, then $\frak{p} = (x_1, x_2, \ldots, x_i)$ for some $i$.  
\end{corollary}


The associated primes of products of monomial primes are known; see \cite[Exercise 3.9]{eisenbud}, but our approach offers an additional benefit:  When $I$ is a principal $Q$-Borel ideal, we can easily read the associated primes and an irredundant primary decomposition from the poset $Q$ and the unique $Q$-generator of $I$. This irredundant primary decomposition is also described in a different framework in the recent paper of Herzog, Rauf, and Vladoiu \cite{HRV} for powers of $I$. 

For $x_i \in Q$, let $A(x_i) = \{ x_j : x_j \leq_Q x_i\}$.  More generally, if $m$ is any monomial, let $A(m) = \{x_j: x_j \leq_Q x_i$ for some $x_i | m\}$.  

\begin{theorem}\label{ass}
Let $I = \Q(m)$ for some monomial $m$.  Then $\frak{p}$ is an associated prime of $S/I$ if and only if $\frak{p} = A(m')$ for some $m'$ dividing $m$ such that $A(m')$ is connected.  
\end{theorem}

\begin{proof}
Write $I = \prod \p^{e_{\p}} = \bigcap \p^{a_{\p}}$. By Proposition~\ref{primefactor}, $e_{\p} > 0$ implies that $\p = A(x_i)$ for some $x_i$ dividing $m$.

By Proposition \ref{qborelprimes}, we need only consider primes which correspond to order ideals of $Q$.  So suppose $A$ is an order ideal of $Q$, and let $\p = A$. Let $T = \{x_k : x_k \mid m, A(x_k) \subseteq \p\}$. If $\p \neq A(m')$ for some $m'$ dividing $m$, we will show that $\p$ is redundant. Put $\q = \sum_{x_k \in T} A(x_k)$, and note that $\q = A(m')$ for some monomial $m'|m$.  Then $\q \subseteq \frak{p}$ and, by Theorem \ref{primeproduct}, $a_{\q} = a_{\p}$, meaning $\frak{p}$ is redundant.  

Now let $m'|m$ and suppose $\p=A(m')$ is not connected. Then there exist monomials $m_1$ and $m_2$ with $m_1 m_2 = m'$ and $A(m_1) \cap A(m_2)  = \emptyset$. Put $\p_1 = A(m_1)$ and $\p_2 = A(m_2)$. By Theorem \ref{primeproduct}, $a_{\p_1} + a_{\p_2} = a_{\p}$ because the exponents $e_{\q}$ are zero for nonprincipal primes; in particular, the exponents are zero when a prime has elements from multiple components.  Because $\frak{p}_1, \frak{p}_2 \subseteq \frak{p}$, it follows that $\frak{p}^{a_{\p}} \supseteq \p_{1}^{a_{\p_1}} \cap \p_{2}^{a_{\p_2}}$ is redundant.  

Finally, let $\frak{p}= A(m')$ for some $m'|m$ such that $A(m')$ is connected. We will construct a monomial $\mu$ such that $\Ann_{S/I}(\mu)=\frak{p}$.  Let $G$ be the graph with vertex set $T=\supp(m')$ and an edge connecting two variables of $T$ if and only if they have a common lower bound in $Q$.  $G$ is connected because $A(m')$ is connected; fix a minimal spanning tree $H$ of $G$.  For each edge $(t_{i},t_{j})$ of $H$, choose a common lower bound $y_{i,j}$ of $t_{i}$ and $t_{j}$.  Set $Y=\prod_{H}y_{i,j}$ and $Z=\prod_{t\in T}t$.  We claim that $\mu=m\frac{Y}{Z}$ is the desired monomial.  Indeed, $\mu x_{i} \in \Q(m)$ if $x_{i}\in A(Z)=A(m')$, and $\mu x_{i}\not\in\Q(m)$ otherwise.
\end{proof}

\begin{remark}
As a consequence of Theorems~\ref{primeproduct} and \ref{ass}, we can read off an irredundant primary decomposition of any principal $Q$-Borel ideal $\Q(m)$ from $m$ and the order ideals of $Q$. Using Theorem~\ref{ass}, we determine which primes are associated. For each of these primes $\frak{p}$, we compute the exponent of $\frak{p}$ in a primary decomposition by summing the exponents $e_i$ on each variable contained in $\frak{p}$. 
\end{remark}

\begin{example}
Let $Q$ be the poset given by the figure below, and let $I = \Q(def)$. 
\begin{figure}[htp]
\centering
\includegraphics[height = 1.1in]{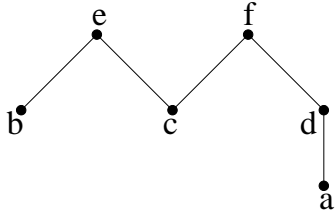}
\caption{The poset $Q$.}\label{path}
\end{figure}

The possible monomial primes associated to $I$ are $A(d) = (a, d)$, $A(e) = (b, c, e)$, $A(f) = A(df) = (a, c, d, f)$, $A(de) = (a, b, c, d, e)$, and $A(ef) = A(def) = (a, b, c, d, e, f)$.  Except for $A(de)$, all these order ideals are connected.  

Thus $\Ass(S/I) = \{ (a, d), (b, c,e), (a, c, d, f), (a, b, c, d, e, f)\}$, and \[ I=(a,d) \cap (b,c,e) \cap (a,c,d,f)^2 \cap (a,b,c,d,e,f)^3 \] is an irredundant primary decomposition of $I$.
\end{example}

\begin{corollary}
Let $I$ be a Borel ideal with a single Borel generator.  That is, $I = \Q(m)$ for some monomial $m$, where $Q = x_1 < x_2 < \ldots < x_n$ is the $n$-element chain.  Then $\frak{p}$ is an associated prime of $I$ if and only if $\frak{p} = (x_1, x_2, \ldots, x_i)$ for some $x_i |m$.  
\end{corollary}

\begin{proof}
Since every order ideal of the chain is connected, Theorem \ref{ass} gives the result.  
\end{proof}

\section{Projective dimension and Cohen-Macaulayness} \label{s.cm}

Throughout this section, let $Q$ be a poset, and let $I=\Q(m) \subset S=k[x_1,\dots,x_n]$ be a principal $Q$-Borel monomial ideal generated by $m$. We assume that all maximal elements of $Q$ divide $m$; if not, then we may make $Q$ smaller by deleting these maximal elements and pass to a smaller ring without changing the minimal generators of $I$.  

We begin by determining the projective dimension of the principal $Q$-Borel ideal $I=\Q(m)$. Recall from Section~\ref{s.preliminaries} that $I$ is polymatroidal (Proposition \ref{polym}). In \cite[Lemma 1.3]{herzogtakayama}, Herzog and Takayama prove that polymatroidal ideals have linear quotients with respect to descending reverse-lex order. For our purposes, it is also convenient to order the monomials in ascending reverse-lex order. (This is equivalent to ordering the monomials in descending graded-lex order with the usual order of the variables reversed, so $x_n > x_{n-1} > \cdots > x_1$.) 

For example, the descending reverse-lex ordering on monomials of degree two in $k[a,b,c]$ is \[ a^2, ab, b^2, ac, bc, c^2, \] while the ascending order starts at $c^2$ and proceeds backwards through the same list.

The next lemma is a result of Francisco and Van Tuyl \cite[Proposition 2.9]{franciscovantuyl}; we include a proof since the argument was omitted in the published version of \cite{franciscovantuyl}.

\begin{lemma}\label{ascendinglinearquo}
Polymatroidal ideals have linear quotients with respect to the ascending reverse-lex order.
\end{lemma}

\begin{proof}
Let $I$ be a polymatroidal ideal, and let $u$ be a minimal generator of $I$. Let $J$ be the ideal generated by all minimal generators $v$ of $I$ less
than $u$ in reverse-lex order. Note that, since $J$ is a monomial ideal, $J:(u)=(\frac{v}{\gcd(u,v)}: v \in J)$. We show that for each $v < u$, there exists
$x_i \in J:(u)$ such that $x_i$ divides $\frac{v}{\gcd(u,v)}$.
                                                                                
Write $u=x_1^{a_1} \cdots x_n^{a_n}$ and $v=x_1^{b_1} \cdots x_n^{b_n}$. Since $v < u$, there exists $i$ such that $a_i < b_i$ but $a_{i+1}=b_{i+1},\dots,$ $a_n=b_n$. We now invoke the dual version of the polymatroidal exchange property (\cite[Lemma 2.1]{herzoghibi}): since $a_i <
b_i$, there exists $j$ with $a_j > b_j$ so that $u'=u \cdot \frac{x_i}{x_j} \in I$. Note that $j < i$, so $u'<u$, and thus $u' \in J$. Because $x_j u'=x_i u$ and $u' \in
J$, we have $x_i \in J:(u)$. The power of $x_i$ in the monomial $\frac{v}{\gcd(u,v)}$ is $b_i-\min(a_i,b_i)=b_i-a_i>0$, and hence $x_i$ divides $\frac{v}{\gcd(u,v)}$.
\end{proof}

\begin{theorem}\label{pdim}
Let $Q$ be a poset and $I=\Q(m)$ for some monomial $m$. Suppose that $Q$ is the maximal poset stabilizing $Q$ and that the maximal elements of $Q$ divide $m$. Then \[ \pd(S/I) = n - \#(\text{connected components of } \, Q) + 1.\]
\end{theorem}

\begin{proof}
By Lemma~\ref{ascendinglinearquo} and \cite[Lemma 1.3]{herzogtakayama}, polymatroidal ideals have linear quotients with respect to either ascending or descending reverse-lex order. Thus we get a minimal free resolution of $I$ by constructing an iterated mapping cone resolution using either order. Let $t$ be the number of connected components of $Q$.

First, we order the minimal monomial generators of $I$ in descending reverse-lex order, so $x_1^d$ is the first monomial in degree $d$. Let $\mu$ be any minimal generator of $I$, and let $I_{\mu}$ be the ideal generated by all minimal generators of $I$ greater than $\mu$ in the reverse-lex order. Let $C$ be any connected component of $Q$, and write $x_{\max(C)}$ for the variable of largest index in $C$ (which divides $m$ by hypothesis).

We claim that for each $C$, $x_{\max(C)} \notin I_{\mu}:\mu$. For any $\nu \in I$ with $\deg \nu = \deg m$, the number of factors of $\nu$ in $C$, counted with multiplicity, is the same as the number of factors of $m$ in $C$, again counted with multiplicity. This follows from the fact that exchanges of variables take place within connected components of $Q$. Suppose $x_{\max(C)} \in I_{\mu}:\mu$. Then $x_{\max(C)} \mu \in I_{\mu}$, and because $I_{\mu}$ is generated in degree one lower, there exists $x_i$ dividing $\mu$ such that $\mu \cdot \frac{x_{\max(C)}}{x_i} \in I_{\mu}$. Because $x_{\max(C)} \in C$, we have $x_i \in C$ as well. Therefore $i < \max(C)$, and $\mu > \mu \cdot \frac{x_{\max(C)}}{x_i}$ in the descending reverse-lex order. But this means that $\mu \cdot \frac{x_{\max(C)}}{x_i} \notin I_{\mu}$, a contradiction. Hence for each connected component $C$ of $Q$ and each minimal generator $\mu$ of $I$, $x_{\max(C)} \notin I_{\mu}:\mu$. Therefore $\pd(S/I) \le n-t+1$.

For the opposite inequality, we order the minimal monomial generators of $I$ in ascending reverse-lex order. Equivalently, we use lex order with $x_n > x_{n-1} > \cdots > x_1$, the usual order of the variables reversed. We let $\mu$ be the minimal generator of $I$ that occurs last in that order. In each connected component $C$, let $x_{\min(C)}$ be the variable of smallest index in $C$. We construct a bipartite graph $G$ as follows. Let $X$ be the set of elements of $C$ dividing $m$ and $Y$ be the set of minimal elements of $C$. Let $X \cup Y$ be the vertex set of $G$. Construct an edge between $x \in X$ and $y \in Y$ if and only if $x$ and $y$ are comparable in $Q$. Without loss of generality, relabel the variables of $C$ (and vertices of $G$) so that for $x_i$ and $x_j$ minimal in $C$, $\dist_G(x_i,x_{\min(C)}) < \dist_G(x_j,x_{\min(C)})$ implies that $i < j$. 


Given a minimal element $x_i$, take a minimum length path \[ x_i \to x_{\ell} \to \cdots \to x_{\min(C)}\] from $x_i$ to $x_{\min(C)}$ in $G$. Let $p$ be the minimal index with $x_p \le_Q x_{\ell}$. It follows from the relabeling that $p < i$. Recall that there is a correspondence between the divisors of $m$ (counted with multiplicity) and the divisors of $\mu$ (counted with multiplicity). By the construction of $\mu$, this correspondence sends every copy of $x_{\ell}$ to a copy of $x_p$. In particular, $x_p$ divides $\mu$. 

Let $I_{\mu}$ be the ideal generated by the minimal monomial generators of $I$ except for $\mu$. Because $\mu \cdot \frac{x_i}{x_p} \in I$ and $i \neq p$, we have $\mu \cdot \frac{x_i}{x_p} \in I_{\mu}$. Thus $\pd(S/I) \ge n-t+1$.
\end{proof}

\begin{proposition}\label{codim}
Let $Q$ be a poset, and let $I=\Q(m)$. Then \[ \codim I = \min_{x_i \mid m} |A(x_i)|,\] where $A(x_i)$ is the order ideal generated by $x_i$.
\end{proposition}

\begin{proof}
The codimension of $I$ is the smallest codimension of an associated prime of $S/I$. By Theorem~\ref{ass}, the associated primes are the connected order ideals of $Q$ generated by nonempty subsets of variables dividing $m$, from which the result is immediate.
\end{proof}

As a consequence, we get a special case of Herzog and Hibi's classification of Cohen-Macaulay polymatroidal ideals \cite[Theorem 4.2]{herzoghibi}.

\begin{corollary}\label{cm}
Let $Q$ be a poset, and suppose $I=\Q(m)$ for some monomial $m$. Assume $Q$ is the maximal poset stabilizing $I$ and that the maximal elements of $Q$ divide $m$. Then $S/I$ is Cohen-Macaulay if and only if one of the following holds:
\begin{enumerate}
\item[(1):] $Q$ is the chain $x_1 <_Q x_2 <_Q \cdots <_Q x_n$, and $m=x_n^{a_n}$ for some $a_n$. (That is, $I$ is a power of the maximal homogeneous ideal.)
\item[(2):] $Q$ is the antichain. (That is, $I$ is a principal ideal.)
\end{enumerate} 
\end{corollary}

Relaxing the assumptions of this section, Corollary \ref{cm} becomes the following:

\begin{theorem}[\cite{herzoghibi}]
Let $Q$ be a poset, and suppose $I=\Q(m)$ is a principal $Q$-Borel ideal. Then $S/I$ is Cohen-Macaulay if and only if one of the following holds:
\begin{enumerate}
\item[(1):] $m=x_{i}^{d}$ for some $i$ and $d$ (and $I$ is a power of the prime ideal $Q(x_{i})$).
\item[(2):] $I$ is a principal ideal.
\end{enumerate} 

\end{theorem}

\begin{proof}[Proof of Corollary~\ref{cm}]
Suppose first that $Q$ is connected. By Theorem~\ref{pdim}, $\pd(S/I) = n$, and thus $S/I$ is Cohen-Macaulay if and only if $\codim I = n$. By Proposition~\ref{codim}, this forces $A(x_i)=\{x_1,\dots,x_n\}$ for all $x_i$ dividing $m$. That is, $m=x_n^{a_n}$ and $x_i <_Q x_n$ for all $i$.
We are assuming that $Q$ is the maximal poset stabilizing $I$, and thus $Q$ is the chain $x_1 <_Q x_2 <_Q \cdots <_Q x_n$. Hence $I=(x_1, \dots, x_n)^{a_n}$.

Now suppose that $Q$ is not connected, and let $t$ be the number of connected components of $Q$, so $\pd(S/I) = n-t+1$. Then $I$ is Cohen-Macaulay if and only if the smallest principal order ideal of $Q$ generated by a divisor of $m$ has cardinality $n-t+1$. In this case, because $Q$ has $t$ connected components and $n$ elements, we have one component with cardinality $\ge n-t+1$, and then the other $t-1$ components must have cardinality one. Therefore $Q$ consists of $t-1$ isolated vertices and a connected component of size $n-t+1$ with a unique maximal element. We are assuming that each maximal element of $Q$ divides $m$. Thus each variable that is isolated in $Q$ divides $m$ and, moreover, the codimension of $I$ is the cardinality of the smallest principal order ideal of $Q$ generated by a divisor of $m$. Thus $n-t+1=\pd(S/I)=\codim I=1$, and so $n=t$. This means that $Q$ consists of $n$ isolated vertices, so $Q$ is the antichain.
\end{proof}


\section{Interpolation between Borel ideals and monomial ideals} \label{s.interpolation}
In this section, we describe procedures for computing free
resolutions and primary decompositions 
of  $\Q$-Borel ideals which specialize to familiar objects in both the
Borel case (i.e., $Q$ is the chain) and the arbitrary monomial case
(i.e., $Q$ is the antichain).  We begin with some remarks on
intersections of $\Q$-Borel ideals.

The intersection of two $\Q$-Borel ideals is again $\Q$-Borel, and can
be computed efficiently in terms of their $\Q$-Borel generators:

\begin{lemma}
Let $I=\Q(T_{1})$ and $J=\Q(T_{2})$, for sets $T_{1}$ and $T_{2}$ of
monomials.  Then $I\cap J = \sum_{t_{1}\in T_{1}, t_{2}\in T_{2}}
\Q(t_{1})\cap \Q(t_{2})$.  
\end{lemma}

In the extremal cases, the intersection of principal (Borel) ideals is
always principal.  This is not the case for general $Q$, as the
following example shows.

\begin{example}\label{noartinian}
Let $S=k[a,b,c]$, with poset structure $a<b$, $a<c$.  Then
for all $k$ we have $\Q(a^{k}b)\cap \Q(a^{k}c)= \Q(a^{k}bc, a^{k+1})$,
and $\Q(a^{k}bc)\cap \Q(a^{k+1})=\Q(a^{k+1}b, a^{k+1}c)$.
\end{example}

In particular, if $I=\Q(T)$, then $\cap_{t\in T}\Q(t)$ need not be a
principal $\Q$-Borel ideal. In fact,
successively intersecting generators as in this example can yield an
infinite decreasing chain of $\Q$-Borel ideals.  Note, however, that
such a chain stabilizes at zero in every degree.

\subsection{Resolutions}
Let $T$ be a set of monomials, and put $I=\Q(T)$.  Recall that $\Q(m)$
has linear resolution for all monomials $m$.  Now let $m\in T$ be a
monomial of minimal degree, and put $I_{m}=\Q(T\smallsetminus m)$.  

The Mayer-Vietoris sequence is
\[
0\to I_{m}\cap \Q(m) \to I_{m}\oplus \Q(m)\to  I\to 0.
\]

Given resolutions for $I_{m}$, $\Q(m)$, and $I_{m}\cap\Q(m)$, the
mapping cone would yield a (usually not minimal) resolution of $I$.
$\Q(m)$ has known linear resolution (see \cite{herzogtakayama}), and $I_{m}$,
having fewer generators than $I$, can in some sense be resolved
inductively.  Unfortunately, as Example \ref{noartinian} shows,
$I_{m}\cap \Q(m)$ is not well-behaved:  An attempt to resolve $I$
inductively in this manner may never terminate.  However, we can
produce a truncated resolution.

\begin{algorithm}\label{resolutionalgorithm}
For a fixed degree $d$, this algorithm produces a complex $\mathbb{F}$
of $R$-modules which satisfies $F_{0}=I$ and $(H_{i}(\mathbb{F}))_{j}=0$ whenever
$j-i\leq d$. 

\begin{description}
\item[Step 1]  Delete any $\Q$-generators of $I$ having degree greater
  than $d$; let $J$ be the resulting ideal.
\item[Step 2]  If $J$ is a principal $\Q$-Borel ideal, then it is
  minimally resolved by the linear quotients on its monomial
  generators (see section \ref{s.preliminaries}).  Otherwise, choose a $\Q$-generator
  $m$ of minimal degree, and let $J_{m}$ be the $\Q$-Borel ideal
  generated by the $\Q$-generators of $J$ other than $m$.  Let
  $\mathbb{G}$ be the minimal resolution of $\Q(m)$.
\item[Step 3]  Observe that $\Hilb_{\leq d}(J_{m})$ and $\Hilb_{\leq
  d}(J_{m}\cap \Q(m))$ are both smaller than $\Hilb_{\leq d}(J)$.
  Inductively apply this algorithm to find complexes $\mathbb{H}$
  and $\mathbb{K}$ that agree with resolutions of $J_{m}$ and
  $\Q(m)\cap J_{m}$ in degrees less than $d$.
\item[Step 4]  Set $\mathbb{F}=\mathbb{G}\oplus
  \mathbb{H}\oplus\mathbb{K}[-1]$.  
\end{description}
\end{algorithm}

\begin{remark}  If $Q$ is the antichain, and $d$ is greater than or
  equal to the degree of the least common multiple of the generators
  of $I$, then Algorithm
  \ref{resolutionalgorithm} produces the Taylor resolution \cite{taylor} of $I$.
\end{remark}

\begin{remark}
The truncated resolution produced by Algorithm
\ref{resolutionalgorithm} is a sum of linear strands corresponding to
monomial generators of the principal $\Q$-Borel ideals arising at
various stages of the algorithm.  If at any stage $(J_{m}\cap \Q(m))$
shares a monomial generator with either $J_{m}$ or $\Q(m)$, we may
obtain a smaller complex by performing a cancellation on the
resulting linear strands.

If $\Q$ is the chain (so $I$ is Borel), $d$ is greater than or equal
to the largest degree of a generator of $I$, and we perform these
cancellations at every opportunity, we get the Eliahou-Kervaire
resolution \cite{EK} of $I$ (see \cite{FMS}).
\end{remark}

\subsection{Irreducible decompositions}
We describe a process for producing an irreducible decomposition of a
$\Q$-Borel ideal $I$.

\begin{definition}
We say that an ideal $J$ is \emph{$Q$-irreducible} if
$J$ has a $Q$-generating set consisting of pure powers of variables.  
\end{definition}

\begin{lemma}  Suppose that $J=\Q(x_{1}^{e_{1}},\dots, x_{n}^{e_{n}})$
  is a $\Q$-irreducible ideal.  (Here, we allow $e_{i}$ to be integers in $[0,\infty]$ with the understanding that $x_{i}^{\infty}=0$.)  Then $J$ has (possibly redundant)
  irreducible decomposition $J=\bigcap (x_{1}^{f_{1}},\dots,
  x_{n}^{f_{n}})$, the intersection taken over exponent vectors
  $(f_{i})$ satisfying the conditions:

\begin{itemize}
\item[(i)] If $e_{k}<\infty$ and $x_{i}\leq_{Q}x_{k}$, then $f_{i}<\infty$.
\item[(ii)] For all $k$, $\sum_{x_{i}\leq_{Q}x_{k}}(f_{i}-1)\leq
  e_{k}-1$.  
\end{itemize}
\end{lemma}
\begin{proof}
Assume without loss of generality that all $e_{i}$ are finite.
Let $T$ be the intersection of irreducible ideals given above.
Suppose $m\in T$, and let $x_{i}$ be a variable dividing $m$.  Then
$\frac{m}{x_{i}}\in T'=\bigcap (x_{1}^{g_{1}},\dots, x_{n}^{g_{n}})$,
where $g_{i}=f_{i}-1$ and $g_{j}=f_{j}$ for $j\neq i$.  Inductively,
$T'=\Q(x_{1}^{a_{1}},\dots, x_{n}^{a_{n}})$, where $a_{j}=e_{j}-1$ whenever
$x_{i}\leq_{Q}x_{j}$ and $e_{j}\geq 1$, and $a_{j}=e_{j}$ otherwise.
We conclude $m\in x_{i}T'=J$.
\end{proof}

Thus to obtain an irreducible decomposition of a $\Q$-Borel ideal $I$,
it suffices to decompose $I$ into $\Q$-irreducible ideals.

\begin{lemma}\label{principalirreducibledecomplemma}
Let $I=Q(m)$ be a principal $Q$-Borel ideal.  Let $z$ be
$Q$-maximal among the variables dividing $m$, and write
$m=z^{e_{z}}\mu\nu$, where every element $x$ of $\supp(\mu)$ satifies
$x\lneqq_{Q}z$, and every element of $\supp(\nu)$ is incomparable to
$z$.  Put $d=\deg(\mu)$.  Then $I=Q(z^{e_{z}+d}\nu)\cap Q(\mu\nu)$.  
\end{lemma}
\begin{proof} Let $J$ denote the intersection, and suppose $f\in J$ is
  a monomial.  To prove $f\in I$, we have to show the existence of a
  matching from the variables dividing $m$ (counted with multiplicity)
  to those dividing $f$, with the property that $x_{i}>_{Q}x_{j}$
  whenever a copy of $x_{i}$ is matched to a copy of $x_{j}$.  We will
  use Hall's marriage theorem.

  Let $X$ be a subset of the divisors of $m$; we will show that, if
  $N(X)$ is the set of of divisors $x_{j}$ of $f$ satisfying
  $x_{j}<_{Q}x_{i}$ for some $x_{i}\in X$, then $|N(X)|\geq
  |X|$.  Write $X=Y\cup Z$, where the elements of $Z$ are all
  comparable to $z$ and the elements of $Y$ are not.  

  First suppose
  $z\in Z$.  Then, since $f\in Q(z^{e_{Z}+d}\nu)$, there is a matching
  from the divisors of this monomial to the divisors of $f$ and we
  have $|X|\leq |Y\cup \{\text{$(e_{z}+d)$ copies of $z$}\}|\leq
  |N(Y\cup \{z\})| = N(X)$.  

  Now suppose that $z\not\in Z$.  Then, since $f\in Q(\mu\nu)$, there
  is a matching from the divisors of this monomial to the divisors of
  $f$ and we have $|X|\leq |N(X)|$ as desired.
\end{proof}

\begin{algorithm}\label{irreducibledecompalg}
This algorithm computes a $\Q$-irreducible decomposition of a
$\Q$-Borel ideal $I=Q(m_{1},\dots, m_{s})$.
\begin{description}
\item[Step 1]  If $I$ is already $\Q$-irreducible, we are done.  If
  not, without loss of generality let $z$ be $Q$-maximal among the
  variables dividing those $Q$-generators of $I$ which are not pure
  powers, and suppose $z$ divides $m_{1}$.
\item[Step 2]  Write $m_{1}=z^{e_{z}}\mu$; set $J=I+Q(z^{e_{z}})$ and
  $K=I+Q(\mu)$.  By Lemma
  \ref{principalirreducibledecomplemma} we have $I=J\cap K$.
\item[Step 3]  Both $J$ and $K$ have a smaller number of
  $Q$-generators which are divisible by $z$ and are not pure powers
  (or a smaller number of variables dividing their non-pure-power
  $Q$-generators), so this algorithm may be recursively applied to
  find $\Q$-irreducible decompositions of $J$ and $K$.
\end{description}
\end{algorithm}

\begin{remark}  When $Q$ is the antichain, Algorithm
  \ref{irreducibledecompalg} is the algorithm for producing
  irreducible decompositions given in \cite[Lemma 5.18]{MS}.
\end{remark}

\begin{remark}  When $Q$ is the chain (so $I$ is Borel), Algorithm
  \ref{irreducibledecompalg} reduces to the following.
\end{remark}

\begin{proposition}
Let $I$ be a Borel ideal, and assume without loss of generality that
$x_{n}$ divides some generators of $I$ but that $I$ contains no power
of $x_{n}$.  Write $I=N+M$, where $N$ is the Borel ideal generated by
the Borel generators of $I$ which are divisible by $x_{n}$, and $M$ is
the Borel ideal generated by those generators which are not.  Let $d$
be the minimal degree of a generator of $N$.  Then
$I=(\Borel(x_{n}^{d})+M)\cap (M+(N:x_{n}^{\infty}))$.  
\end{proposition}

\section{$Y$-Borel ideals} \label{s.y-borel}

In our final section, we investigate ideals which are Borel with respect to the poset $Y$ defined by the relations $x_1 <_Y \cdots <_Y x_t <_Y y$ and $x_t <_Y z$. Let $I$ be a $Y$-Borel ideal in $S=k[x_1,\dots,x_t,y,z]$. The minimal generators of $I$ not divisible by $z$ form a Borel ideal in the ring $k[x_1,\dots,x_t,y]$. Because $Y$ is similar to the chain on all the variables of $S$, it is natural to believe that the minimal free resolution of $I$ shares many properties of an Eliahou-Kervaire resolution. We determine the minimal free resolution of $I$ in this section, providing evidence for our broader belief that if a poset $Q$ is close to the chain, a $Q$-Borel ideal will behave much like a Borel ideal.

We first recall the Eliahou-Kervaire resolution.

\begin{notation}[\cite{PS}]  Given a monomial $m\in k[x_{1},\dots, x_{n}]$, set $\max(m)=\max\{i:x_{i} \text{ divides } m\}$ and $\min(m) =\min\{i: x_{i}\text{ divides }m\}$.  
Let $I$ be a Borel ideal of $k[x_{1},\dots, x_{t}]$.  Then, if $\mu$ is a monomial of $I$, there exists a unique factorization $\mu=\mu_{1}\mu_{2}$ such that $\mu_{1}$ is a generator of $I$ and $\max(\mu_{1})\leq \min(\mu_{2})$.  We say that $\mu_{1}$ and $\mu_{2}$ are the \emph{beginning} and \emph{end} of $\mu$, respectively, and write $\Beg(\mu)=\mu_{1}$ and $\END(\mu)=\mu_{2}$.
\end{notation}

\begin{theorem}[\cite{EK}, \cite{PS}]
Let $I$ be a Borel ideal of $k[x_{1},\dots, x_{t}]$.  Then the minimal free resolution of $I$ has basis given by the ``Eliahou-Kervaire symbols'' $[m,\alpha]$, where $m$ is a monomial generator of $I$ and $\alpha$ is a squarefree monomial with $\max(\alpha)\lneqq \max(m)$.  The symbol $[m,\alpha]$ has homological degree $\deg(\alpha)$ and multidegree $m\alpha$.  The differential is 
\[
d([m,\alpha])= \sum_{i=1}^{\deg \alpha} (-1)^{1+i}\alpha_{i}[m,\frac{\alpha}{\alpha_{i}}] - \sum_{i=1}^{\deg\alpha}(-1)^{1+i} 
\END (m\alpha_{i}) [\Beg (m\alpha_{i}),\frac{\alpha}{\alpha_{i}}],
\]
where $\alpha_{i}$ is the $i^{\text{th}}$ variable (in lex order) dividing $\alpha$, and the symbol $[\Beg (m\alpha_{i}),\frac{\alpha}{\alpha_{i}}]$ is treated as zero if it is not an Eliahou-Kervaire symbol (i.e., if $\max(\frac{\alpha}{\alpha_{i}})\geq\max(\Beg (m\alpha_{i}))$).  
\end{theorem}

The minimal free resolution of a $Y$-Borel ideal is similar, but nonlinear syzygies between powers of $y$ and $z$ appear because the variables are incomparable in in $Y$.

\begin{theorem}\label{t:yborel}
Let $I$ be a $Y$-Borel ideal.  Then the minimal free resolution of $I$ has basis given by the symbols $[m,\alpha]$ and $[m,\alpha y^{k_{m}}]$, where $m$ is a monomial generator of $I$, $\alpha$ is a squarefree monomial of $k[x_{1},\dots, x_{t}]$ satisfying $\max(\alpha)\lneqq \max(m)$, and $k_{m}$ is minimal such that $\frac{m}{z}y^{k_{m}}\in I$.  The symbol $[m,\alpha]$ has homological degree $\deg(\alpha)$ and  multidegree $m\alpha$.  The symbol $[m,\alpha y^{k_{m}}]$ has homological degree $1+\deg(\alpha)$ and multidegree $m\alpha y^{k_{m}}$.  The differential of the symbol $[m,\alpha]$ is exactly the Eliahou-Kervaire differential.  The differential of the symbol $[m,\alpha y^{k_{m}}]$ is 
\begin{align*}
d([m,\alpha y^{k_{m}}])=&\sum_{i=1}^{\deg \alpha} (-1)^{1+i}\alpha_{i}[m,\frac{\alpha}{\alpha_{i}}y^{k_{m}}] + (-1)^{\deg(\alpha)}y^{m_{k}}[m,\alpha]\\
&-\Bigg(\sum_{i=1}^{\deg\alpha}(-1)^{1+i} 
\END (m\alpha_{i}) [\Beg (m\alpha_{i}),\frac{\alpha}{\alpha_{i}}y^{k_{m}}]\\
&\quad\quad + (-1)^{\deg(\alpha)}\END (my^{k_{m}}) [\Beg (my^{k_{m}},\alpha]\Bigg),
\end{align*}
where, as in the Eliahou-Kervaire resolution, symbols that don't exist are treated as zero.
\end{theorem}

The proof of Theorem \ref{t:yborel} is by induction on the largest power of $z$ appearing in a generator of $I$.  We need some new notation and a lemma.

\begin{notation}
Given a $Y$-Borel ideal $I$, put $I_{z}=\frac{1}{z}(I\cap (z))$, and let $I_{1}$ be the ideal generated by the monomial generators of $I$ which are not divisible by $z$.  Because $I$ is an ideal, we have $I_{1}\subset I_{z}$.  Note also that $I=I_{1}+zI_{z}$, that $I_{1}$ is Borel in $k[x_{1},\dots, x_{t},y]$, that $I_{z}$ is $Y$-Borel, and that Theorem \ref{t:yborel} applies to $I_{z}$ by induction.
\end{notation}

\begin{lemma}\label{smallend} Suppose that $m$ is a monomial generator of $I_{1}$.  Because $m\in I_{z}$, we may calculate $\Beg(m)$ and $\END(m)$ in $I_{z}$.  
If $x_{\max(m)}=x_{i}$ for some $i\leq t$, then $\END(m)=1$ or $\END(m)=x_{\max(m)}$.  If $x_{\max(m)}=y$, then $\END(m)=1$ or $\END(m)=y^{k_{mz}}$.
\end{lemma}
\begin{proof} Suppose that $x_{\max(m)}=x_{i}$ for some some $i\leq t$, and that $x_{j}x_{i}$ divides $\END(m)$.  Then $\frac{mz}{x_{j}x_{i}}\in I$.  Applying the Borel move sending $z$ to $x_{j}$, we have $\frac{m}{x_{i}}\in I$.  Since this monomial is not divisible by $z$, it is in $I_{1}$, contradicting the assumption that $m$ was a minimal generator for $I_{1}$.  The proof in the other case is similar.
\end{proof}

\begin{proof}[Proof of Theorem \ref{t:yborel}]
We induct on the largest power of $z$ appearing in a generator of $I$.  In the base case, $I\subset k[x_{1},\dots, x_{t},y]$ is Borel.
In general, recall the Mayer-Vietoris sequence:
\[
0\to I_{1}\cap zI_{z} \to I_{1}\oplus zI_{z} \to I_{1}+zI_{z} \to 0.
\]
This simplifies to 
\[
0\to zI_{1} \to I_{1}\oplus zI_{z} \to I \to 0,
\]
which is isomorphic to 
\[
0 \to I_{1}(z^{-1})\stackrel{\begin{pmatrix}z\\ -1\end{pmatrix}}{\longrightarrow} I_{1}\oplus I_{z}(z^{-1}) \stackrel{\begin{pmatrix}1 &z\end{pmatrix}}{\longrightarrow}I\to 0.
\]
$I$ is resolved (nonminimally) by the mapping cone of this short exact sequence.  We will analyze the cancellations in this mapping cone, allowing us to reduce to the minimal resolution.  First, we introduce some extra notation.

By induction, a basis for the resolution of $I_{z}$ is given by symbols of the form $[m,\alpha]$ and $[m,\alpha y^{k_{m}}]$ for generators $m$ of $I_{z}$.  Above $I_{z}(z^{-1})$ in the mapping cone, we will refer to these symbols as $[zm,\alpha]$ and $[zm, \alpha y^{k_{m}}]$.  

$I_{1}$ is Borel in $k[x_{1},\dots, x_{t},y]$, so a basis for its resolution is given by the Eliahou-Kervaire symbols $[m,\alpha]$.  Above $I_{1}$ in the middle term of the Mayer-Vietoris sequence, we refer to this by its natural name $[m,\alpha]$.
Above $I_{1}(z^{-1})$, we will refer to the symbol $[m,\alpha]$ as $[z\Beg (m),\alpha\END(m)]$, where beginning and end are computed in the ideal $\frac{1}{z}I_{z}$. However, by Lemma~\ref{smallend}, $\END(m)$ is divisible by at most one variable.

Since $I_{1}(z^{-1})$ has componentwise linear resolution, the only opportunity for cancellation in the mapping cone is in the lift of the map from $I_{1}(z^{-1})$ to $I_{z}(z^{-1})$.  The lift of this map sends the symbol $[z\Beg(m), \alpha\END(m)]$ to $(-1)^{1+\deg \alpha}\END(m)[z\Beg(m),\alpha]$ if $\END(m)$ is not divisible by $y$, and to $(-1)^{1+\deg \alpha}\frac{\END(m)}{y^{k_{mz}}}[z\Beg(m),\alpha y^{k_{mz}}]$ if $\END(m)$ is divisible by $y^{k_{mz}}$.  In the first case, there is cancellation if and only if $\END(m)=1$, i.e., if and only if $m$ is a minimal generator for $I_{z}$.  In the second case, there is cancellation if and only if $\END(m)=y^{k_{mz}}$, i.e., if and only if $\frac{m}{y^{k_{mz}}}$ is a minimal generator for $I_{z}$.

The uncancelled symbols fall into three categories:
\begin{enumerate}
\item[(1):]  Symbols of the form $[m,\alpha]$ where $z$ does not divide $m$, arising from the resolution of $I_{1}$ in the middle term.

\item[(2):]  Symbols of the form $[mz, \alpha]$ or $[mz, \alpha y^{k_{mz}}]$ where $m$ is a minimal generator for $I_{z}$ but not for $I_{1}$, arising from the resolution of $I_{z}$ in the middle term.  

\item[(3):]Symbols of the form $[z\Beg(m), \alpha\END(m)]$ where $z$ does not divide $m$, and $m$ is a minimal generator for $I_{1}$ but not for $I_{z}$, arising from the resolution of $I_{1}$ in the intersection term.
\end{enumerate}

Observe that these categories are disjoint.  Thus it suffices to show that if $\mu$ is a monomial generator of $I$, and $\beta$ is a squarefree monomial with $\max(\beta)\lneqq \max(\mu)$, then $[\mu,\beta]$ and $[\mu,\beta y^{k_{\mu}}]$ each fall in one of the three categories above.  

If $\mu$ is not divisible by $z$, then set $m=\mu$ and $\alpha=\beta$, and observe that $[\mu, \beta] = [m,\alpha]$ in the first category.  In this case, $[\mu,\beta y^{k_{\mu}}]$ is nonsensical.

If $\mu$ is divisible by $z$, set $\nu=\frac{\mu}{z}$.  If $\max(\beta)\lneqq \max(\nu)$, set $m=\nu$ and $\alpha=\beta$ and observe that $[\mu,\beta]$ and $[\mu,\beta y^{k_{\mu}}]$ are $[mz, \alpha]$ or $[mz, \alpha y^{k_{mz}}]$ from the second category.  If not, set $m=\nu x_{\max(\beta)}$ and $\alpha = \frac{\beta}{x_{\max(\beta}}$, and observe that $[\mu,\beta]$ and $[\mu,\beta y^{k_{\mu}}]$ are $[z\Beg(m), \alpha\END(m)]$ from the third category. \end{proof}

\bigskip

\noindent \textbf{Acknowledgements}: Markus Vasquez made substantial contributions at the beginning of our efforts to investigate these problems, and we thank him for his work. We also thank Craig Huneke and Sonja Mapes for helpful conversations. This work was partially supported by grants from the Simons Foundation (\#199124 to Christopher Francisco and \#202115 to Jeffrey Mermin).  



\begin{thebibliography}{9}

\bibitem[BS]{bayerstillman} D. Bayer and M. Stillman, A theorem on refining division orders by the reverse lexicographic order. \emph{Duke Math. J.} {\bf 55} (1987), no. 2, 321--328.

\bibitem[CH]{concaherzog} A. Conca and J. Herzog, Castelnuovo-Mumford regularity of products of ideals. \emph{Collect. Math.} {\bf 54} (2003), no. 2, 137--152.

\bibitem[E]{eisenbud} D. Eisenbud, \emph{Commutative Algebra with a View Toward Algebraic Geometry}, Springer, New York, NY, 1995.

\bibitem[EK]{EK} S. Eliahou and M. Kervaire, Minimal resolutions of some monomial ideals. \emph{J. Algebra} {\bf 129} (1990), no. 1, 1–-25.

\bibitem[FMS]{FMS} C. A. Francisco, J. Mermin, and J. Schweig, Borel generators. \emph{J. Algebra} {\bf 332} (2011), 522--542.

\bibitem[FVT]{franciscovantuyl} C. A. Francisco and A. Van Tuyl, Some families of componentwise linear monomial ideals. \emph{Nagoya Math. J.} {\bf 187} (2007), 115--156.

\bibitem[HH]{herzoghibi} J. Herzog and T. Hibi, Cohen-Macaulay polymatroidal ideals. \emph{European J. Combin.} {\bf 27} (2006), no. 4, 513--517.

\bibitem[HRV]{HRV} J. Herzog, A. Rauf, and M. Vladoiu, The stable set of associated prime ideals of a polymatroidal ideal. Preprint, 2011. {\tt arXiv:1109.5834}.

\bibitem[HT]{herzogtakayama} J. Herzog and Y. Takayama, Resolutions by mapping cones. The Roos Festschrift volume, 2. \emph{Homology Homotopy Appl.} {\bf 4} (2002), no. 2, part 2, 277--294.

\bibitem[MS]{MS} E. Miller, B. Sturmfels, {\it Combinatorial Commutative Algebra.} GTM 227, Springer-Verlag, New York, 2004.

\bibitem[PS]{PS} I. Peeva and M. Stillman, The minimal free resolution of a Borel ideal. \emph{Expo. Math.} {\bf 26} (2008), no. 3, 237–-247. 

\bibitem[S]{ec1} R. Stanley, \emph{Enumerative Combinatorics Vol 1}, Cambridge University Press, Cambridge, 1997.  

\bibitem[T]{taylor} D. Taylor, \emph{Ideals generated by monomials in an R-sequence}, Ph.D. Thesis, University of Chicago, 1960. 

\end{thebibliography}
\end{document}